\documentclass[10pt, twoside, a4paper]{article}
\let\finishall\relax\let\Finishall\relax\let\getprepared\relax
\ifx\TestIngCommand\undefined\relax
\else\input ../../../proc/ppreamb-book.tex 
  \input init.tex \fi
\let\TestIngCommand\undefined

\newtheorem{remark}{Remark}

\newtheorem{conj}{Conjecture}
\usepackage{amsmath,amscd,amssymb}                                         
\usepackage{graphics}                                                     
\newtheorem{theo}{Theorem}                                                 
\newtheorem{lem}{Lemma}                                                    
\newtheorem{cor}{Corollary}                                                
\newtheorem{defi}{Definition}
\newskip\ttglue\ttglue=.5em plus.25em minus.15em                           
\def\firstname#1{\def\FIRSTNAME{#1}\ignorespaces}
\def\lastname#1{\def\LASTNAME{#1}\ignorespaces}
\def\middleinitial#1{\def\MIDDLEINI{#1}\ignorespaces}
\def\department#1{\def\DEPARTMENT{#1}\ignorespaces}
\def\institute#1{\def\INSTITUTE{#1}\ignorespaces}
\def\address#1{\def\ADDRESS{#1}\ignorespaces}
\def\country#1{\def\COUNTRY{#1}\ignorespaces}
\def\otheraffiliation#1{\def\OTHERAFFILIATION{#1}\ignorespaces}
\def\email#1{\def\EMAIL{#1}\ignorespaces}
\newcount\autcount\autcount=0                                              
\newcount\affcount\affcount=0                                              
\newcount\numcount\newcount\nummcount                                      
\newcount\nummmcount\newcount\nummmmcount                                  
\def\writename#1#2{\ \kern-1ex\hbox{
  \csname AUthor\the#1\endcsname\                                          
  \edef\TESTSTR{}\expandafter\ifx\csname auTHor\the#1\endcsname\TESTSTR    
  \else\csname auTHor\the#1\endcsname.\ \fi                                
  \csname authOR\the#1\endcsname$^{\csname AFF\the#1\endcsname}$
  \expandafter\ifx\csname corr\number#1\endcsname\relax                    
  \else\thanks{Corresponding author.}\ \fi                                 
  }\ifnum#1<#2, \else\ \kern-1ex\fi}
\def\writeemail#1{
  \nummcount=0\relax\nummmcount=0\relax                                    
  \loop\ifnum\nummcount<\autcount\advance\nummcount by1\relax              
    {\expandafter\ifnum\csname AFF\the\nummcount\endcsname=#1\relax        
    \global\advance\nummmcount by1\fi}\repeat                              
  \nummcount=0\relax\nummmmcount=0\relax                                   
  \loop\ifnum\nummcount<\autcount\advance\nummcount by1\relax              
    {\expandafter\ifnum\csname AFF\the\nummcount\endcsname=#1\relax        
    \global\advance\nummmmcount by1\relax\def\blank{}\expandafter          
    \ifx\csname EMAIL\the\nummcount\endcsname\blank(no e-mail)
    \else\csname EMAIL\the\nummcount\endcsname                             
    \fi                                                                    
    \ifnum\nummmmcount<\nummmcount; \fi\fi}\repeat}
\long\def\BeginAuthorList#1\EndAuthorList{#1\relax                         
  \author{\vbox{\hsize=390pt\noindent\numcount=0\relax                     
    \loop\ifnum\numcount<\autcount\advance\numcount by1\relax              
      \writename{\numcount}{\autcount}
      \repeat}\\[2mm]                                                      
    \vbox{\small\numcount=0\relax                                          
      \loop\ifnum\numcount<\affcount\advance\numcount by1\relax            
        \vbox{{\count0=\numcount\relax                                     
          \loop\expandafter\ifnum\csname AFF\the\count0\endcsname
            <\numcount\relax\advance\count0 by1\relax\repeat               
          $^{\csname AFF\the\count0\endcsname}$}
        \def\BLANK{}\expandafter\ifx\csname DEPT\the\numcount\endcsname    
          \BLANK                                                           
          \else\csname DEPT\the\numcount\endcsname, \fi                    
        \csname INST\the\numcount\endcsname,                               
        \csname ADDR\the\numcount\endcsname,                               
        \csname COUN\the\numcount\endcsname                                
        \edef\TEST{}\expandafter\ifx\csname OTHE\the\numcount\endcsname
          \TEST                                                            
          .\else;\break\csname OTHE\the\numcount\endcsname.\fi}
        \vbox{\writeemail{\numcount}}
        \repeat}\\}}
\expandafter\def\csname x1\endcsname{}
\expandafter\def\csname x2\endcsname{}
\expandafter\def\csname x3\endcsname{}
\expandafter\def\csname x4\endcsname{}
\expandafter\def\csname x5\endcsname{}
\expandafter\def\csname x6\endcsname{}
\expandafter\def\csname x7\endcsname{}
\expandafter\def\csname x8\endcsname{}
\expandafter\def\csname x9\endcsname{}
\def\Author#1#2{\global\advance\autcount by1\relax#2                       
  \expandafter\edef\csname AUthor\the\autcount\endcsname{\FIRSTNAME}
  \expandafter\edef\csname auTHor\the\autcount\endcsname{\MIDDLEINI}
  \expandafter\edef\csname authOR\the\autcount\endcsname{\LASTNAME}
  \expandafter\edef\csname EMAIL\the\autcount\endcsname{\EMAIL}
  \let\tempera\"\def\"{\string\"}\expandafter\ifx\csname x\DEPARTMENT
    \endcsname\relax                                                       
    \global\advance\affcount by1\relax\let\"\tempera                       
    \expandafter\edef\csname DEPT\the\affcount\endcsname{\DEPARTMENT}
    \expandafter\edef\csname INST\the\affcount\endcsname{\INSTITUTE}
    \expandafter\edef\csname ADDR\the\affcount\endcsname{\ADDRESS}
    \expandafter\edef\csname COUN\the\affcount\endcsname{\COUNTRY}
    \expandafter\edef\csname OTHE\the\affcount\endcsname{\OTHERAFFILIATION}
    \expandafter\edef\csname AFF\the\autcount\endcsname{\the\affcount}
  \else\expandafter\edef\csname AFF\the\autcount\endcsname{\DEPARTMENT}
  \fi\let\"\tempera\ignorespaces}
\def\CorrespondingAuthor#1#2{
  \expandafter\xdef\csname corr\number#1\endcsname{cor}
  \Author#1{#2}}
\def\PaperTitle#1{\title{\bf#1}}
\def\Category#1{\ignorespaces}
\def\keywords#1{{\noindent \emph{Keywords:}                                
  \def\BLANK{}\def\TEST{#1}\ifx\BLANK\TEST(n/a).\else#1\fi}}
\getprepared                                                               
\begin{document}                                                           

\PaperTitle{Distribution of Small Values of Bohr Almost Periodic Functions with Bounded Spectrum}
\Category{(Pure) Mathematics}

\date{}

\BeginAuthorList
  \Author1{
    \firstname{Wayne}
    \lastname{Lawton}
    \middleinitial{M}   
    \department{Department of the Theory of Functions, Institute of Mathematics and Computer Science}
    \institute{Siberian Federal University}
    \address{Krasnoyarsk}
    \country{Russian Federation}
    \otheraffiliation{}
    \email{wlawton50@gmail.com}}
\EndAuthorList
\maketitle
\thispagestyle{empty}
\begin{abstract}
If $f$ is a nonzero Bohr almost periodic function on $\mathbb R$ with a bounded spectrum we prove there exist $C_f > 0$ and integer $n > 0$ such that for every $u > 0$ the mean measure of the set $\{\, x \, : \, |f(x)| < u \, \}$ is less than $C_f\, u^{1/n}.$ For trigonometric polynomials with $\leq n + 1$ frequencies we show that $C_f$ can be chosen to depend only on $n$ and the modulus of the largest coefficient of $f.$ We show this bound implies that the Mahler measure $M(h),$ of the lift $h$ of $f$ to a compactification $G$ of $\mathbb R,$ is positive and discuss the relationship of Mahler measure to the Riemann Hypothesis.
\end{abstract}
%

\finishall
\section{Distribution of Small Values}
\label{statement}
$\mathbb{N} \mathrel{\mathop :}= \{1,2,...\}, \mathbb{Z}, \mathbb{R}, \mathbb{C},
\mathbb{T} \mathrel{\mathop :}= \{\, z \in \mathbb C \, : \, |z| = 1 \ \}$
are the natural, integer, real, complex and circle group numbers, $\textrm C_b(\mathbb R)$ is the $\textrm C^*$-algebra of bounded continuous functions and
$\chi_\omega : \mathbb R \rightarrow \mathbb T, \ \omega \in \mathbb R$ the homomorphisms
$\chi_{\, \omega}(x) \mathrel{\mathop :}= e^{i \omega x},\ \omega \in \mathbb R.$ A finite sum $f = \sum_{\, \omega} a_{\, \omega} \, \chi_\omega$ with distinct $\omega$ is called a trigonometric polynomial with height
$H_f \mathrel{\mathop :}= \max_{\, \omega} |a_{\,\omega}|$
and they comprise the algebra $\textrm T(\mathbb R)$ of trigonometric polynomials. Bohr \cite{bohr1} defined the $\textrm C^*$-algebra $\textrm U(\mathbb R)$ of uniformly almost periodic functions to be the closure of $\textrm{T}(\mathbb R)$ in $\textrm C_b(\mathbb R)$ and proved that their means
$m(f) \, \mathrel{\mathop :}= \lim_{L \rightarrow \infty} (2L)^{-1} \int_{-L}^L f(t) dt$
exist. The Fourier transform $\widehat f \, : \, \mathbb R \rightarrow \mathbb C$ of $f \in \textrm U(\mathbb R)$ is
$\widehat f(\omega) \mathrel{\mathop :}= m(f \, \overline \chi_\omega \, )$
and its spectrum
$\Omega(f) \mathrel{\mathop :}= \hbox{ support } \widehat f.$
If $f$ is nonzero then $\Omega(f)$ is nonempty and countable and we say
$f$ has bounded spectrum if its bandwidth $b(f) \in [0,\infty],$ defined by
$b(f) \mathrel{\mathop :}= \sup \Omega(f) - \inf \Omega(f),$
is finite.
We observe that if $S \subseteq \mathbb R$ is defined by a finite number of inequalities involving functions in $\textrm U(\mathbb R)$ then
$m(S) \, \mathrel{\mathop :}= \lim_{L \rightarrow \infty} (2L)^{-1}
\hbox{measure } [-L,L] \cap S$
exists and define $J_f : (0,\infty) \rightarrow [0,1]$ by
\begin{equation}\label{Jf}
    J_f(u) \mathrel{\mathop :}= m \left(\, \{ \, x \in \mathbb R \, : \, |f(x)| < u \, \}\, \right)
\end{equation}
\begin{theo}\label{thm1}
If $f \in {\textrm U}(\mathbb R)$ is nonzero and has a bounded spectrum then there exist $C_f > 0$ and $n \in \mathbb N$ such that:
\begin{equation}\label{ineq1}
    J_f(u) \leq C_f \, u^{\frac{1}{n}}, \ \ u > 0.
\end{equation}
There exists a sequence $C_n$ such that if $f \in \textrm T(\mathbb R)$ has $n+1$ frequencies then
\begin{equation}\label{ineq2}
        J_f(u) \leq C_n \, H_f^{-\frac{1}{n}} \, u^{\frac{1}{n}}, \ \ u > 0.
\end{equation}
\end{theo}
{\bf Proof} For $f \in \textrm U(\mathbb R), \omega \in \mathbb R, k \in \mathbb N, u > 0$ define $\Xi_{f,\omega, k,u}\, ,  K_f \, : \, (0,\infty) \rightarrow [0,1]$ by
\begin{equation}\label{Xi}
    \Xi_{f,\omega, k,u}(v) \mathrel{\mathop :}=
    m \{ \, x \in \mathbb R \, : \, |f(x)| < u, \,
    |(\, \chi_{\omega}f)^{(j)}(x)\, | < v^j, j = 1,...,k \, \},
\end{equation}
\begin{equation}\label{Kf}
    K_f(u) \mathrel{\mathop :}=
    \inf_{\omega \in \mathbb R} \, \inf_{k \in \mathbb N} \, \inf_{v > 0} \,
    \left[ 3\, {\sqrt 2}\, \pi^{-1}\, b(f)\, k \, v^{-1}\, u^{\frac{1}{k}} \ + \ \Xi_{f,\omega, k,u}(v)\right].
\end{equation}
We first prove Theorem 1 assuming the following result which we prove latter.
\begin{lem}\label{lem1}
Every nonzero $f \in \textrm U(\mathbb R)$ with bounded spectrum satisfies
$J_f \leq K_f.$
\end{lem}
We observe that for every $\omega \in \mathbb R$ and every $a \in \mathbb R \backslash \{0\},$ if $h(x) = \chi_{\, \omega}(x)\, f(ax)$
then $J_h = J_f$ and $K_h = K_f.$ Without loss of generality we can assume that $\Omega(f) \subset [-\frac{b(f)}{2},\frac{b(f)}{2}].$ If $b(f) = 0$ then $f = c$ and $J_f(u) \leq |c|^{-1}u.$ If $\Omega(f) > 0$ then
Bohr \cite{bohr2} proved that $f$ extends to an entire function $F$ of exponential type $\frac{b(f)}{2},$ and Boas \cite{boas1}, (\cite{boas2}, p. 11, Equation 2.2.12) proved that
\begin{equation}\label{Boaslimit}
    \limsup_{k \rightarrow \infty} |f^{(k)}(x)|^{\frac{1}{k}} = \frac{b(f)}{2}
\end{equation}
uniformly in $x.$ Therefore for any $v_0 > \frac{b(f)}{2}$ there exists $k \in \mathbb N$ such that $\Xi_{f,0,k,u}(v_0) = 0$ so Lemma \ref{lem1} implies $J_f$ satisfies (\ref{ineq1}) with $C_f = 3\, {\sqrt 2}\, \pi^{-1}\, b(f)\, k\, v_0^{-1}$ and $n = k.$ This proves the first assertion. To prove the second we assume, without loss of generality, that $b(f) = 1,$ $\Omega(f) \subset [0,1]$ and
$$f(x) = \sum_{j = 1}^{n+1} a_j \, e^{i \, \omega_j \, x}, \ \
0 = \omega_1 < \cdots < \omega_{n+1} = 1, \ \
H_f = \max \{\, |a_j| \, : \, j = 2,...,n+1 \, \}.$$
Define $C_1 \mathrel{\mathop :}= \frac{1}{2}.$  If $n = 1$ and $f$ has $n + 1 = 2$ terms
and $f = a_0 + a_1\chi_1$ with $|a_1| = H_f$ and $h = H_f(1 - \chi_1),$
then
$
    J_f(u) \leq J_{h}(u) = (2/\pi) \sin^{-1}(u/(2H_f)) \leq C_1 H_f^{-1} u
$
therefore (\ref{ineq2}) holds for $n = 1.$ For $n \geq 2$ we assume by induction that (\ref{ineq2}) holds for $n-1$ and therefore, since $f^{(1)}$ has $n$ terms and $H_{f^{(1)}} = H_f,$ it follows that for all $v > 0,$
\begin{equation}\label{Jfbound}
J_{f^{(1)}}(v) \leq C_{n-1}\, H_f^{\frac{1}{n-1}}\, v^{\frac{1}{n-1}},
\end{equation}
\begin{equation}\label{Xibound}
\Xi_{f,0,1,u}(v) \leq C_{n-1}\,  H_f^{\frac{1}{n-1}} \, v^{\frac{1}{n-1}}.
\end{equation}
Therefore Lemma 1 with $\omega = 0,\, b(f) = k = 1$ gives
\begin{equation}\label{Jubound}
    J_f(u) \leq \inf_{v > 0} \left[3\, {\sqrt 2}\, \pi^{-1}\,  \, v^{-1}\, u +
    C_{n-1} H_f^{\frac{1}{n-1}} \, v^{\frac{1}{n-1}}\right] =
    C_n \, H_f^{\frac{1}{n}}\, u^{\frac{1}{n}}
\end{equation}
\begin{equation}\label{Cn}
 \hbox{where  }   C_n \mathrel{\mathop :}= C_{n-1}^{1-\frac{1}{n}}\, [3\sqrt 2 \, \pi^{-1}\, (n-1)]^{\frac{1}{n}} \, n(n-1)^{-1}.
\end{equation}
\begin{remark}
\label{200million}
Computation of 200 million terms shows that $n^{-1}C_n \rightarrow 0.900316322$
\end{remark}
\begin{conj}\label{conj1}
    In (\ref{ineq2}) $C_n$ can be replaced by a bounded sequence.
\end{conj}
\begin{lem}\label{lem2}
If $\phi : [a,b] \rightarrow \mathbb{C}$ is differentiable
and $\phi^{\prime}([a,b])$ is contained in a quadrant then
\begin{equation}
\label{muI}
    b-a \leq 2{\sqrt 2} \ \frac{\max |\phi|([a,b])}{\min |\phi^{\prime}|([a,b])}
\end{equation}
\end{lem}
{\bf Proof of Lemma 2} We first proved this result in (\cite{lawton4}, Lemma 1) where we used it to give a proof, of a conjecture of Boyd \cite{boyd} about monic polynomials related to Lehmer's conjecture \cite{lehmer}, which was reviewed in (\cite{everestward}, Section 3.5) and extended to monic trigonometric polynomials in (\cite{lawton9}, Lemma 2).
The triangle inequality
$|\phi^{\prime}| \leq |\Re\, \phi^{\prime}| + |\Im\, \phi^{\prime}|$
gives
$$(b-a)\min |\phi^{\prime}|([a,b]) \leq \int_{a}^{b} |\phi^{\prime}(y)| \, dy \leq \int_{a}^{b}
\left( \, |\Re \, \phi^{\prime}(y)| +  |\Im \, \phi^{\prime}(y)| \, \right) \, dy.$$
Since $\phi^{\prime}([a,b])$ is contained in a quadrant of $\mathbb C$ there exist
$c, d \in \{1,-1\}$ such that $|\Re \, \phi^{\prime}(y)| = c\, \Re \, \phi^{\prime}(y)$
and $|\Im \, \phi^{\prime}(y)| = d \, \Im \, \phi^{\prime}(y)$ for all $y \in [a,b].$
Therefore
$$\int_{a}^{b}
\left( \, |\Re \, \phi^{\prime}(y)| +  |\Im \, \phi^{\prime}(y)| \, \right) \, dy
= (c\, \Re \phi(b) + d\, \Im \phi(b)) - (c\, \Re \phi(a) + d\, \Im  \phi(a)).$$
The result follows since the right side
is bounded above by $2{\sqrt 2} \, \max |\phi|([a,b]).$
\\ \\
{\bf Proof of Lemma 1}
Assume that $f \in \textrm U(\mathbb R)$ is nonzero. We may assume without loss of generality that $\Omega(f) \subset [-\frac{b(f)}{2},\frac{b(f)}{2}].$ For
$k \in \mathbb N, u > 0, v > 0$ we define the set
\begin{equation}\label{S}
    S_{f,k,u,v} \mathrel{\mathop :}= \{ \, x \in \mathbb R \, : \, |f(u)| < u, \, \max_{j \in \{1,...,k\}} \, |f^{(j)}(x)|^{\frac{1}{j}} \geq v \, \}.
\end{equation}
We observe that the set of functions in $\textrm U(\mathbb R)$ whose spectrums are in $[-\frac{b(f)}{2},\frac{b(f)}{2}]$ is closed under differentiation, and define
$
    s(f,k,u,v) \mathrel{\mathop :}= m(S_{f,k,u,v}).
$
\begin{equation}\label{s2}
 \hbox{It suffices to prove that}\ \ \ \    s(f,k,u,v) \leq 3\, {\sqrt 2}\, \pi^{-1}\, b(f)\, k \, v^{-1}\, u^{\frac{1}{k}}.
\end{equation}
Define
$\gamma_j \mathrel{\mathop :}= u^{\frac{k-j}{k}} \, v^j, j \in \{0,...,k\},$
and $\mathcal{I} \mathrel{\mathop :}=$ set of closed intervals $I$ satisfying, for some $j \in \{0,1,...,k-1\},$ the following three properties:
\begin{enumerate}
\item $f^{(j+1)}(I)$ is a subset of a closed quadrant,
\item $\max |f^{(j)}|(I) \leq \gamma_j$ and
$\min |f^{(j+1)}|(I) \geq \gamma_{j+1},$
\item $I$ is maximum with respect to properties 1 and 2.
\end{enumerate}
Define $\mathcal{E} \mathrel{\mathop :}= $ set of endpoints of intervals in $\mathcal{I}$, and
\begin{equation}\label{psi}
\psi \mathrel{\mathop :}= \prod_{j = 0}^{k-1}(\Re f^{(j+1)})(\Im f^{(j+1)})(|f^{(j)}(x)|^2 - \gamma_j^2)(|f^{(j+1)}(x)|^2 - \gamma_{j+1}^2).
\end{equation}
\begin{equation}
 \hbox{Lemma \ref{lem2} implies that}\ \ \ \    \hbox{length }(I) \leq 2{\sqrt 2}\, \frac{\gamma_k}{\gamma_{k+1}} = 2{\sqrt 2}\, v^{-1}\, u^{\frac{1}{k}}, \ \ I \in \mathcal{I},
\end{equation}
\begin{equation}\label{smallcontained}
 \hbox{and (\ref{S}) and Property 3 implies that}\ \ \ \    S_{f,k,u,v} \subset \bigcup_{I \in \mathcal{I}} I.
\end{equation}
Clearly $\psi = \Psi|_{\mathbb R}$ where $\Psi$ is the product of $6k$ entire functions each having bandwidth $b(f)$ so a theorem of Titchmarsh \cite{titchmarsh1} implies that the density of real zeros of $\Psi$ is bounded above by
$3\pi^{-1}\, b(f)\, k.$
Property 3 implies that all points in $\mathcal{E}$ are zeros of $\Psi$
so the upper density of intervals in $\mathcal{I}$ is bounded by
$\frac{3}{2}\, \pi^{-1}\, b(f)\, k.$
Combining these facts gives
$s(f,k,u,v) \leq (\frac{3}{2}\pi^{-1}\, b(f)\, k)\, (2 \sqrt 2 v^{-1}u^{\frac{1}{k}}) =
3 \sqrt 2 \pi^{-1}\, b(k)\, k\, v^{-1} u^{\frac{1}{k}}$
which proves (\ref{s2}) and concludes the proof of Lemma 1.
\\ \\
For $p \in [1,\infty)$ Besicovitch \cite{besicovitch} proved that the completion $\textrm B^p(\mathbb R)$ of $\textrm U(\mathbb R)$ with norm $(m(|f|^p))^{\frac{1}{p}}$ is a subset of $L_{loc}^p(\mathbb R).$
For $x \geq 0$ we define $\log^{+}(x) \mathrel{\mathop :}= \log (\max \{1,x\}) \in [0,\infty),$ $\log^{-}(x) \mathrel{\mathop :}= \log (\min \{1,x\}) \in [-\infty,0],$ and $|x|_j \mathrel{\mathop :}= \max \{|x|, \frac{1}{j}\}$ for
$j \in \mathbb N.$
\begin{cor}\label{Bp}
    If $f \in \textrm U(\mathbb R)$ satisfies (\ref{ineq1}), then
    $\log^- \circ \, |f| \in \textrm B^p(\mathbb R),$
    \begin{equation}\label{lowbound}
        m(\, |\log^- \circ \, |f|\, |^p\, ) \leq \int_0^1 |\log(u)|^p \, d \, C_f \, u^{\frac{1}{n}} = C_f\, n^p\, \Gamma(p),
    \end{equation}
    and $\log \circ |f| \in B^p(\mathbb R).$
\end{cor}
{\bf Proof of Corollary 1} Since the means of the functions $\log^{-} \circ |f|_j\, |^p$ are nondecreasing and bounded by the right side of (\ref{lowbound}), the sequence $\log^{-} \circ |f|_j $ is a Cauchy sequence in $B^p(\mathbb R)$ so it converges to a function $\eta \in B^p(\mathbb R).$ Therefore $\log^{-} \circ |f| = \eta$ since it is the pointwise limit of $\log^{-} \circ |f|_j$ and $\eta \in L_{loc}^p(\mathbb R).$
The last fact follows since $\log = \log^+ + \log^{-}.$
\section{Compactifications and Hardy Spaces}
\begin{defi}\label{defi1}
A compactification of $\mathbb R$ is a pair $(G,\theta)$ where $G$ is a compact abelian group and $\theta \, : \, \mathbb R \rightarrow G$ is a continuous homomorphism with a dense image.
\end{defi}
$\textrm C(G)$ is the set of continuous functions on $G$ and $L^p(G), p \in [1,\infty)$ are Banach spaces. If $h \in \textrm C(G)$ then $f \mathrel{\mathop :}= h \circ \theta \in \textrm U(\mathbb R)$ since by a theorem of Bochner \cite{bochner} every sequence of translates of $f$ has a subsequence that converges uniformly. We call $h$ the lift of $f$ to $G.$
The Pontryagin dual \cite{pontryagin} $\widehat G$ of a compact abelian group $G$ is the discrete group of continuous homomorphisms $\chi \, : \, G \rightarrow \mathbb T$ under pointwise multiplication. Bohr proved the existence of a compactification $(\mathbb{B},\theta)$ such that $\textrm U(\mathbb R) = \{ \, h \circ \theta\, : \, h \in \textrm C(\mathbb{B}) \, \}.$ The group $\mathbb{B}$ is nonseparable and $\widehat {\mathbb B}$ is isomorphic to $\mathbb R_{\, d} \mathrel{\mathop :}=$ real numbers with the discrete topology.
\begin{lem}\label{lem3}
For every $f \in \textrm U(\mathbb{R})$ there exists a compactification $(G(f),\theta)$ and $h \in \textrm C(\mathbb{R})$
such that $f = h \circ \theta.$ The group $G(f)$ is separable.
\end{lem}
{\bf Proof of Lemma 3} If $f \in \textrm U(\mathbb R)$ is nonzero its spectrum $\Omega(f)$ is nonempty and countable so the product group $\mathbb T^{\, \Omega(f)}$ is compact and separable. The function
$\theta \, : \, \mathbb R \rightarrow \mathbb T^{\, \Omega(f)}$ defined by $\theta(x)(\omega) \mathrel{\mathop :}= \chi_{\omega}(x)$ is a continuous homomorphism. Define $G(f) \mathrel{\mathop :}= \overline {\theta(\mathbb R)}.$ Then $(G(f),\theta)$ is a compactification. The function $\widetilde h \, : \, \theta(\mathbb R) \rightarrow \mathbb C$ defined by
$\widetilde h(\theta(x) \mathrel{\mathop :}= f(x)$ is uniformly continuous so extends to a unique function $h \, : \, G \rightarrow \mathbb C$ and $f = h\circ \theta.$
\begin{lem}\label{lem4}
If $(G,\theta)$ is a compactification, $h \in \textrm C(G),$ $f = h \circ \theta,$ and $\log \circ |f| \in B^p(\mathbb R),$ then $\log \circ |h| \in L^p(G)$ and $\int_{G} |\log \circ |h| \, |^p = m(|\log \circ |f| \, |^p).$
\end{lem}
{\bf Proof of Lemma 4} The theorem of averages (\cite{arnold}, p. 286) implies that
\begin{equation}\label{dominated}
\int_G |\log^{-} \circ |h|_j\, |^p =
m(|\log^{-} \circ |f|_j\, |^p) \leq m(|\log^{-} \circ |f| |^p).
\end{equation}
The result follows from Lebesgue's monotone convergence theorem since the sequence $|\log \circ |h|_j\, |^p$ is nondecreasing, converges pointwise to $|\log \circ |h|\, |^p$ pointwise and by (\ref{dominated}) their integrals are
uniformly bounded.
\begin{defi}\label{Fourier}
The Fourier transform $\mathfrak{F} \, : \, L^1(G) \rightarrow \ell^\infty(\widehat G)$ is defined by
$\mathfrak{F}(h)(\chi) \mathrel{\mathop :}= \int_G f \, \overline \chi.$
\end{defi}
We define the spectrum $\Omega(h) \mathrel{\mathop :}= \hbox{support }\mathfrak{F}(h).$
The Hausdorff-Young theorem \cite{hausdorff, young} implies that the restrictions give bounded operators $\mathfrak{F} \, : \, L^p(G) \rightarrow \ell^q(\widehat G)$ for $p \in [1,\infty)$ and $p^{-1} + q^{-1} = 1.$
\begin{defi}\label{Order}
A compactification $(G,\theta)$ induces an injective homomorphism $\xi \, : \widehat G \rightarrow \mathbb R,$ $\xi(\chi) \mathrel{\mathop :}= \omega$ where $\chi \circ \theta = \chi_\omega,$ by which we will identity $\widehat G$ as a subset of $\mathbb R$ with the same archimedian order. Therefore if $h \in \textrm C(G)$ is the lift of $f \in \textrm U(\mathbb R),$ then $\Omega(h) = \Omega(f).$ The compactification gives Hardy spaces $H^p(G,\theta) \mathrel{\mathop :}= \{ \, h \in L^p(G) \, : \, \Omega(h) \subset [0,\infty)\, \}, p \in [1,\infty].$
\end{defi}
\begin{defi}\label{outer}
A function $h \in H^p(G,\theta)$ is outer if
$\int_G h \neq 0,$ $\log \circ |h| \in L^1(G),$  and
\begin{equation}\label{outerint}
    \int_G \log \circ  |h| = \log \, \left| \, \int_G h \, \right|.
\end{equation}
A function $h \in H^p(G,\theta)$ is inner if $|h| = 1.$
\end{defi}
A polynomial $h$ is outer iff it has no zeros in the open unit disk since
a formula of Jensen \cite{jensen} gives
$
\int_G \log \circ \, |h| =  \log |h(0)| - \sum_{h(\lambda) = 0}
    \log^{-} (|\lambda|).
$
Beurling \cite{beurling} proved that a function $h \in H^2(\mathbb{T})$ admits a factorization $h = h_o\, h_i,$ with $h_o$ outer and $h_i$ inner, iff $\log \circ \, |h| \in L^1(\mathbb{T}).$
\\ \\
Let $(G,\theta)$ be a compactification. If $h \in \textrm C(G)$ has a bounded spectrum $\Omega(h) \subset [0,\infty)$ then $f = h \circ \theta$ extends to an entire function $F$ bounded in the upper half plane. We observe that if $F$ has no zeros in the upper half plane, then $\chi_{-b(h)/2}F$ is the Ahiezer spectral factor \cite{ahiezer} of the entire function $F(z)\overline {F(\overline z)}.$
\begin{conj}\label{conj2}
$h$ above is outer iff $F$ has no zeros in the open upper half plane.
\end{conj}
\section{Mahler Measure and the Riemann Hypothesis}
\begin{defi}\label{defi2}
For $G$ is a compact abelian group the Mahler measure \cite{mahler1, mahler2}
of $h \in L^1(G)$ is $M(h) \mathrel{\mathop :}= \exp \left(\, \int_G \log \circ |h|\, \right) \in [0,\infty).$ We also define
$M^{\pm}(h) \mathrel{\mathop :}= \exp \left(\, \int_G \log^{\pm} \circ |h|\, \right).$
\end{defi}
Since $M(h) = M^+(h)M^-(h)$
and $M^+(h) \in [\, 1,\, \max \{1,||h||_{\infty}\}\, ],$ it follows that
$M(h) > 0$ iff $\log^{-} \circ |h| \in L^1(G)$ and then
$M^{-}(h) = \exp \left(-||\, \log^{-} \circ |h|\, ||_1\, \right).$ Lemma \ref{lem4} implies that this condition holds whenever $h \in \textrm C(G)$ is nonzero and $\Omega(h)$ is bounded.
\begin{defi}\label{PhiN}
    For $N \in \mathbb N, \, \Phi_N \mathrel{\mathop :}= $ product of the first $N$ cyclotomic polynomials.
\end{defi}
Amoroso (\cite{amoroso}, Theorem 1.3) proved that the Riemann Hypothesis is equivalent to
\begin{equation}\label{RH1}
    \log M^+(\Phi_N) \ll_\epsilon N^{\frac{1}{2} + \epsilon},  \ \  \epsilon > 0.
\end{equation}
Define $f_N  \mathrel{\mathop :}= \Phi_N \circ \chi_1 \in \textrm U(\mathbb R)$
and define $J_{f_N} : (0,\infty) \rightarrow [0,1]$ by (\ref{Jf}).
Jensen's formula implies that $M(\Phi_N) = 1$ therefore
\begin{equation}\label{RH2}
    \log M^+(\Phi_N) = - \int_0^1 \log(u) \, d \, J_{f_N}(u).
\end{equation}
The bounds that we obtained for $J_{f}$ in (\ref{ineq1}) and (\ref{ineq2}) were exceptionally crude and totally inadequate to obtain (\ref{RH1}). When deriving (\ref{ineq2}) for general polynomials we used the bound
(\ref{Xibound})
$\Xi_{f,0,1,u}(v) = m(\{\, x\, :\, |f(x)| < u, |f^{(1)}(x)| < v\, \} \leq
m(\{\, x\, :\, |f^{(1)}(x)| < v\, \}.$ Conjecture (\ref{conj1}) was based on our intuition that a smaller upper bound holds. We suspect that much smaller upper bounds hold for specific sequences of polynomials as illustrated by the following examples. Construct sequences of height $1$ polynomials
\begin{equation}\label{PQ}
P_n(z) \mathrel{\mathop :}= 1 + z + \cdots + z^n \ ; \ Q_n(z) \mathrel{\mathop :}= \binom{n}{\, [n/2]\, }^{-1}(1+z)^n
\end{equation}
and $p_n \mathrel{\mathop :}= P_n \circ \chi_1, q_n \mathrel{\mathop :}= Q_n \circ \chi_1.$ Both polynomials have maxima at $z = 1,$
$||P_n||_{\infty} = n+1,$ Stirling's approximation gives $||Q_n||_{\infty} \approx
\sqrt {\pi \, n/2}$ for large $n,$ and for $u \in (0,1]$
\begin{equation}\label{Jpn}
J_{p_n}(u) \leq \frac{2}{\pi}\sin^{-1}\left(\, \min \{\, 1, u\, \} \, \right) \leq u \Rightarrow \log (M^{-}(P_n)) > -1,
\end{equation}
\begin{equation}\label{Jqn}
J_{q_n}(u) = \frac{2}{\pi} \sin^{-1}\left(
\min \left\{
\, 1,  \frac{1}{2} \binom{n}{\, [n/2]\, }^{\frac{1}{n}} \, u^{\frac{1}{n}}\,
\right\}
\right)
\geq \frac{2}{\pi} u^{\frac{1}{n}} \Rightarrow \log (M^{-}(Q_n)) < -\frac{2n}{\pi}.
\end{equation}
Differences between these polynomials arise from their root discrepancy. Those of $P_n$ are nearly evenly spaced. Those of $Q_n,$ all at $z = -1,$ have maximally discrepancy.
\begin{conj}\label{conj3}
    If $R_n$ is a polynomial with $n+1$ terms and height $H(R_n) = 1$ then
    $M^-(P_n) \leq M^-(R_n) \leq M^-(Q_n).$
\end{conj}
The roots of $\Phi_N$ have the form $\exp (2\pi i a_k), k = 1,...,\deg \Phi_N$ where $a_k$ are the Farey series consisting of rational numbers in $[0,1)$ whose denominators are $\leq N.$ Bounds on the discrepancy of the Farey series were shown by Franel \cite{franel} and by Landau \cite{landau} to imply the Riemann Hypothesis.
The relationship between the discrepancy of roots of a polynomial and its coefficients, and the distributions of roots of entire functions have been extensively studied since the seminal paper by
Erd\"{os} and Tur\'{a}n \cite{erdos} and the extensive work by Levin and his school \cite{levin}. We suggest that investigation of the functions
$\Xi_{f,\omega,k,v,u}$ in (\ref{Xi}) and derived functions $K_f$ in (\ref{Kf}) functions may further elucidate how the distribution of small values of polynomials and entire functions depend on their coefficients and roots.
\\ \\
{\bf Acknowledgment} The author thanks Professor August Tsikh
for insightful discussions.

\Finishall
\end{document}